\input amstex
\documentstyle {amsppt}
\magnification=\magstephalf \rightheadtext {Metric fixed point
theory} \topmatter
\title Two results in metric fixed point theory
 \endtitle
\author  Daniel Reem, Simeon Reich and Alexander J. Zaslavski \endauthor
\address Department of Mathematics, The Technion-Israel Institute of
Technology, 32000 Haifa, Israel \endaddress \email
dream\@tx.technion.ac.il; sreich\@tx.technion.ac.il; \newline
ajzasl\@tx.technion.ac.il \endemail \abstract We establish two
fixed point theorems for certain mappings of contractive type. The
first result is concerned with the case where such mappings take a
nonempty, closed subset of a complete metric space $X$ into $X$,
and the second with an application of the continuation method to
the case where they satisfy the Leray-Schauder boundary condition
in Banach spaces.
\endabstract
\endtopmatter
\document

\head 1. Introduction  \endhead

In spite of its simplicity (or perhaps because of it), the Banach
fixed point theorem still seems to be the most important result in
metric fixed point theory. As far as we know, the first
significant generalization of Banach's theorem was obtained in
1962 by E. Racotch [5] who replaced Banach's strict contraction
with contractive mappings, that is, with those mappings which
satisfy condition (1.1) below. Such mappings, as well as many
modifications, were studied and used by many authors. Recently, a
renewed interest in contractive mappings has arisen. See, for
example, [6, 7] where genericity and well-posedness results were
established. Another important topic in fixed point theory is the
search for fixed points of nonself-mappings. In the present paper
we combine these two themes by proving two fixed point theorems
for contractive nonself-mappings. In Theorem 1 we present a new
sufficient condition for the existence and approximation of the
unique fixed point of a contractive mapping which maps a nonempty,
closed subset of a complete metric space $X$ into $X$. In Theorem
2 we present an application of the continuation method to
contractive mappings which satisfy the Leray-Schauder boundary
condition on a subset of a Banach space with a nonempty interior.

\proclaim {Theorem 1} Let $K$ be a nonempty, closed subset of a
complete metric space $(X,\rho)$. Assume that $T:K \to X$
satisfies
$$\rho(Tx,Ty) \le \phi(\rho(x,y))\rho(x,y) \text { for each } x,y
\in K, \tag1.1$$ where $\phi:[0,\infty) \to [0,1]$ is a
monotonically decreasing function such that $\phi(t)<1$ for all
$t>0$.

Assume that $K_0 \subset K$ is a nonempty, bounded set with the
following property:

(P1) For each natural number $n$, there exists $x_n \in K_0$ such
that $T^ix_n$ is defined for all $i=1,\dots,n$.

Then

(A) the mapping $T$ has a unique fixed point $\bar x $ in $K$;

(B) For each $M,\epsilon
>0$, there exist $\delta >0$ and a natural number $k$
such that for each integer $n \ge k$ and each sequence
$\{x_i\}_{i=0}^n \subset K$ satisfying
$$\rho(x_0,\bar x)\le M  \text { and }  \rho(x_{i+1},Tx_i) \le \delta,\;
i=0,\dots,n-1,$$ we have
$$\rho(x_i,\bar x)\le \epsilon, \; i=k,\dots,n.
\tag1.2$$
\endproclaim

Let $G$ be a nonempty subset of a Banach space $(Y,||\cdot||)$. In
[3] J. A. Gatica and W. A. Kirk proved that if $T:\overline{G}\to
Y$ is a strict contraction, then $T$ must have a unique fixed
point $x_1$, under the additional assumptions that the origin is
in the interior Int$(G)$ of $G$ and that $T$ satisfies a certain
boundary condition known as the Leray-Schauder condition:
$$Tx\not = \lambda x \quad \forall x\in
\partial G,\quad \forall\lambda>1. \tag{L-S}$$
Here $G$ is not necessarily convex or bounded. Their proof was
nonconstructive. Later, M. Frigon, A. Granas and Z. E. A. Guennoun
[2], and M. Frigon [1] proved that if $x_t$ is the unique fixed
point of $tT$, then, in fact, the mapping $t\to x_t$ is Lipschitz,
so it gives a partial way to approximate $x_1$. Our second result
extends these theorems to the case where $T$ merely satisfies
(1.1).

\proclaim{Theorem 2} Let $G$ be a nonempty subset of a Banach
space $Y$ with $0\in Int(G)$. Suppose that $T:\overline{G}\to X$
is nonexpansive and that it satisfies condition (L-S). Then for
each $t\in [0,1)$, the mapping $tT:\overline{G}\to X$ has a unique
fixed point $x_t\in Int(G)$ and the mapping $t\to x_t $ is
Lipschitz on $[0,b]$ for any $0<b<1$. If, in addition, $T$
satisfies (1.1), then it has a unique fixed point $x_1\in
\overline{G}$ and the mapping $t\to x_t$ is continuous on $[0,1]$.
In particular, $x_1=\lim_{t\to 1^{-}}x_t$.
\endproclaim

Our paper is organized as follows. Part (A) of Theorem 1 is proved
in Section 2 while Section 3 is devoted to the proof of Theorem
1(B). Finally, we prove Theorem 2 in Section 4.

\head 2. Proof of Theorem 1(A) \endhead

 The uniqueness of $\bar x$ is obvious. To establish its
existence let $x_n \in K_0$ be, for each natural number $n$, the
point provided by property (P1). Fix $\theta_0 \in K$. Since $K_0$
is bounded, there is $c_0>0$ such that
$$\rho(\theta,z) \le c_0 \text { for all } z \in K_0. \tag2.1$$
Let $\epsilon >0$. We will show that there exists a natural number
$k$ such that the following property holds:

(P2) If $n > k$ is an integer and if an integer $i$ satisfies $k
\le i<n$, then
$$\rho(T^ix_n,T^{i+1}x_n) \le \epsilon. \tag2.2$$

Let us assume the contrary. Then for each natural number $k$,
there exist natural numbers $n_k$ and $i_k$ such that
$$k \le i_k<n_k \text { and }
 \rho(T^{i_k}x_{n_k},T^{i_k+1}x_{n_k})>\epsilon.
\tag2.3$$ Choose a natural number $k$ such that
$$k>(\epsilon(1-\phi(\epsilon)))^{-1}(2c_0+\rho(\theta,T\theta)).
\tag2.4$$  By (2.3) and (1.1),
$$\rho(T^ix_{n_k},T^{i+1}x_{n_k})>\epsilon,\; i=0,\dots,i_k.
\tag2.5$$ (Here we use the notation that $T^0z=z$ for al $z \in
K$.) It follows from (1.1), (2.5) and the monotonicity of $\phi$
that for all $i=0,\dots,i_k-1$,
$$\rho(T^{i+2}x_{n_k},T^{i+1}x_{n_k}) \le
\phi(\rho(T^{i+1}x_{n_k},T^ix_{n_k}))\rho(T^{i+1}x_{n_k},T^ix_{n_k})$$
$$\le \phi(\epsilon)\rho(T^{i+1}x_{n_k},T^ix_{n_k})$$
and
$$\rho(T^{i+2}x_{n_k},T^{i+1}x_{n_k})-\rho(T^{i+1}x_{n_k},T^ix_{n_k})$$
$$\le
(\phi(\epsilon)-1)\rho(T^{i+1}x_{n_k},T^ix_{n_k})<-(1-\phi(\epsilon))\epsilon.
\tag2.6$$ Inequalities (2.6) and (2.3) imply that
$$-\rho(x_{n_k},Tx_{n_k}) \le
\rho(T^{i_k+1}x_{n_k},T^{i_k}x_{n_k})-\rho(x_{n_k},Tx_{n_k})=
$$
$$=\sum_{i=0}^{i_k-1}[\rho(T^{i+2}x_{n_k},T^{i+1}x_{n_k})-
\rho(T^{i+1}x_{n_k},T^ix_{n_k})] \le -(1-\phi(\epsilon)\epsilon)
i_k \le-k(1-\phi(\epsilon)\epsilon
$$
and
$$k(1-\phi(\epsilon))\epsilon \le \rho(x_{n_k},Tx_{n_k}). \tag2.7$$
In view of (2.7), (1.1) and (2.1),
$$k(1-\phi(\epsilon))\epsilon \le \rho(x_{n_k},Tx_{n_k})$$
$$ \le
\rho(x_{n_k},\theta)+\rho(\theta,T\theta)+\rho(T\theta,Tx_{n_k})
\le c_0+\rho(\theta,T\theta)+c_0$$ and
$$k \le
(\epsilon(1-\phi(\epsilon)))^{-1}(2c_0+\rho(\theta,T\theta)).$$
This contradicts (2.4). The contradiction we have reached proves
that for each $\epsilon >0$, there exists a natural number $k$
such that (P2) holds.

Now let $\delta >0$. We show that there exists a  natural number
$k$ such that the following property holds:

(P3) If $n > k$ is an integer and if integers $i,j$ satisfy $k \le
i,j <n$, then
$$\rho(T^ix_n,T^jx_n) \le \delta.$$

To this end, choose a positive number
$$\epsilon <4^{-1}\delta(1-\phi(\delta)). \tag2.8$$
We have already shown that there exists a natural number $k$ such
that (P2) holds.

Assume that natural numbers $n,i$ and $j$ satisfy
$$n > k \text { and } k\le i,j<n. \tag2.9$$
We claim that $\rho(T^ix_n,T^jx_n) \le \delta.$

Assume the contrary. Then
$$\rho(T^ix_n,T^jx_n)>\delta. \tag2.10$$
By (P2), (2.9), (1.1), (2.10) and the monotonicity of $\phi$,
$$\rho(T^ix_n,T^jx_n) \le
\rho(T^ix_n,T^{i+1}x_n)+\rho(T^{i+1}x_n,T^{j+1}x_n)+
\rho(T^{j+1}x_n,T^jx_n)$$ $$ \le
\epsilon+\rho(T^{i+1}x_n,T^{j+1}x_n)+\epsilon \le
2\epsilon+\phi(\rho(T^ix_n,T^jx_n))\rho(T^ix_n,T^jx_n)$$ $$ \le
2\epsilon +\phi(\delta)\rho(T^ix_n,T^jx_n).$$ Together with (2.8)
this implies that
$$\rho(T^ix_n,T^jx_n) \le 2\epsilon(1-\phi(\delta))^{-1}<\delta,$$
a contradiction. Thus we have shown that for each $\delta
>0$, there exists a natural number $k$ such that (P3) holds.

Let $\epsilon >0$. We will show that there exists a natural number
$k$ such that the following property holds:

(P4) If the integers $n_1,n_2 \ge k$, then
$\rho(T^kx_{n_1},T^kx_{n_2}) \le \epsilon$.

Choose a natural number $k$ such that
$$k>((1-\phi(\epsilon))(\epsilon))^{-1}4c_0 \tag2.11$$
and assume that the integers $n_1$ and $n_2$ satisfy $$n_1,n_2 \ge
k. \tag2.12$$ We claim that $\rho(T^kx_{n_1},T^kx_{n_2}) \le
\epsilon$. Assume the contrary. Then
$$\rho(T^kx_{n_1},T^kx_{n_2})>\epsilon.$$
Together with (1.1) this implies that
$$\rho(T^ix_{n_1},T^ix_{n_2})>\epsilon,\; i=0,\dots,k. \tag2.13$$
By (1.1), (2.13) and the monotonicity of $\phi$ we have for
$i=0,\dots,k-1$,
$$\rho(T^{i+1}x_{n_1},T^{i+1}x_{n_2}) \le
\phi(\rho(T^ix_{n_1},T^ix_{n_2}))\rho(T^ix_{n_1},T^ix_{n_2})$$
$$\le \phi(\epsilon)\rho(T^ix_{n_1},T^ix_{n_2})$$
and
$$\rho(T^{i+1}x_{n_1},T^{i+1}x_{n_2})-\rho(T^ix_{n_1},T^ix_{n_2})
\le (\phi(\epsilon)-1)\rho(T^ix_{n_1},T^ix_{n_2}) \le
-(1-\phi(\epsilon))\epsilon.$$ This implies that
$$-\rho(x_{n_1},x_{n_2}) \le
\rho(T^kx_{n_1},T^kx_{n_2})-\rho(x_{n_1},x_{n_2})$$
$$=\sum_{i=0}^{k-1}[\rho(T^{i+1}x_{n_1},T^{i+1}x_{n_2})-\rho(T^ix_{n_1},T^ix_{n_2})]
\le -k(1-\phi(\epsilon))\epsilon.$$ Together with (2.1) this
implies that
$$k(1-\phi(\epsilon))\epsilon \le \rho(x_{n_1},x_{n_2}) \le
\rho(x_{n_1},\theta)+\rho(\theta,x_{n_2})\le 2c_0.$$ This
contradicts (2.11). Thus we have shown that
$$\rho(T^kx_{n_1},T^kx_{n_2})  \le \epsilon.$$
In other words, there exists a natural number $k$ for which (P4)
holds.

Let $\epsilon >0$. By (P4), there exists a natural number $k_1$
such that
$$\rho(T^{k_1}x_{n_1},T^{k_1}x_{n_2}) \le \epsilon/4 \text  { for
all integers } n_1,n_2 \ge k_1. \tag2.14$$ By (P3), there exists a
natural number $k_2$ such that $$\rho(T^ix_n,T^jx_n) \le
\epsilon/4 \text { for all natural numbers } n,j,i \text {
satisfying } k_2 \le i,j< n. \tag2.15$$

Assume now that the natural numbers $n_1,n_2, i$ and $j$ satisfy
$$n_1,n_2 > k_1+k_2,\; i,j \ge k_1+k_2,\; i<n_1,\; j <
n_2. \tag2.16$$ We claim that
$$\rho(T^ix_{n_1},T^jx_{n_2}) \le
\epsilon.$$ By (2.14), (2.16) and (1.1)
$$\rho(T^{k_1+k_2}x_{n_1},T^{k_1+k_2}x_{n_2}) \le
\rho(T^{k_1}x_{n_1},T^{k_1}x_{n_2}) \le \epsilon/4. \tag2.17$$ In
view of (2.16) and (2.15),
$$\rho(T^{k_1+k_2}x_{n_1},T^ix_{n_1}) \le \epsilon/4 \text { and }
\rho(T^{k_1+k_2}x_{n_2},T^jx_{n_2})\le \epsilon/4.$$ Together with
(2.17) these inequalities imply that
$$\rho(T^ix_{n_1},T^jx_{n_2}) \le
\rho(T^ix_{n_1},T^{k_1+k_2}x_{n_1})$$ $$+
\rho(T^{k_1+k_2}x_{n_1},T^{k_1+k_2}x_{n_2})+
\rho(T^{k_1+k_2}x_{n_2},T^jx_{n_2}) <\epsilon.$$ Thus we have
shown that the following property holds:

(P5) For each $\epsilon
>0$, there exists a natural number $k(\epsilon)$ such that
$$\rho(T^ix_{n_1},T^jx_{n_2})\le \epsilon $$ for all natural numbers
$n_1,n_2 \ge k(\epsilon)$, $i \in [k(\epsilon),n_1)$ and $j \in
[k(\epsilon),n_2).$

Consider the two sequences $\{T^{n-2}x_n\}_{n=2}^{\infty}$ and
$\{T^{n-1}x_n\}_{n=2}^{\infty}$. Property (P5) implies that both
of them are Cauchy and that
$$\lim_{n
\to\infty}\rho(T^{n-1}x_n,T^{n-2}x_n)=0.$$ Therefore there exists
$\bar x \in K$ such that
$$\lim_{n \to \infty} \rho(\bar x,T^{n-2}x_n)=\lim_{n \to
\infty}\rho(\bar x,T^{n-1}x_n)=0.$$ Since the mapping $T$ is
continuous, $T\bar x=\bar x$ and assertion (A) is proved.

\head 3. Proof of Theorem 1(B) \endhead

For each $x \in X$ and $r>0$, set
$$B(x,r)=\{y \in X:\; \rho(x,y) \le r\}. \tag3.1$$
Choose $\delta_0>0$ such that
$$\delta_0<M(1-\phi(M/2))/4.\tag3.2$$
Assume that
$$y \in K\cap B(\bar x,M),\; z \in X \text { and } \rho(z,Ty)\le \delta_0.
\tag3.3$$ By (3.3) and (1.1),
$$\rho(\bar x,z)\le \rho(\bar x,Ty)+\rho(Ty,z)\le \rho(T\bar
x,Ty)+\delta_0\le \phi(\rho(\bar x,y))\rho(\bar x,y)+\delta_0.
\tag3.4$$ There are two cases:
$$\rho(y,\bar x) \le M/2; \tag3.5$$
$$\rho(y,\bar x)>M/2. \tag3.6$$
Assume that (3.5) holds. By (3.4), (3.5) and (3.2),
$$\rho(\bar x,z) \le \rho(\bar x,y)+\delta_0\le M/2+\delta_0<M.
\tag3.7$$ If (3.6) holds, then by (3.4), (3.3), (3.6) and the
monotonicity of $\phi$,
$$\rho(\bar x,z)\le \delta_0+\phi(M/2)\rho(\bar x,y)\le
\delta_0+\phi(M/2)M$$ $$<(M/4)(1-\phi(M/2))+\phi(M/2)M\le M.$$
Thus $\rho(\bar x,z) \le M$ in both cases.

We have shown that
$$\rho(\bar x,z)\le  M \text { for each } z \in X \text { and } y
\in K \cap B(\bar x,M) \text { satisfying } \rho(z,Ty)\le
\delta_0.\tag3.8$$ Since $M$ is any positive number, we conclude
that there is $\delta_1>0$ such that
$$\rho(\bar x,z)\le \epsilon \text { for each } z \in X \text {
and } y \in K \cap B(\bar x,\epsilon) \text { satisfying }
\rho(z,Ty)\le \delta_1. \tag3.9$$ Choose a positive number $\delta
$ such that
$$\delta <\min\{\delta_0,\delta_1,\;
\epsilon(1-\phi(\epsilon))4^{-1}\} \tag3.10$$ and a natural number
$k$ such that
$$k>4(M+1)(1-\phi(\epsilon)\epsilon)^{-1}+4. \tag3.11$$
Let $n \ge k$ be a natural number and assume that $\{x_i\}_{i=0}^n
\subset K$ satisfies
$$\rho(x_0,\bar x) \le M  \text { and } \rho(x_{i+1},Tx_i) \le
\delta,\;i=0,\dots,n-1. \tag3.12$$ We claim that (1.2) holds. By
(3.8), (3.12) and the inequality $\delta <\delta_0$ (see (3.10)),
$$\{x_i\}_{i=0}^k \subset B(\bar x,M). \tag3.13$$
Assume that (1.2) does not hold. Then there is an integer $j$ such
that
$$j \in \{k,n\} \text { and } \rho(x_j,\bar x)>\epsilon. \tag3.14$$
By (3.14), (3.12), (3.9) and (3.10),
$$\rho(x_i,\bar x)>\epsilon,\;
i=0,\dots,j. \tag3.15$$ Let $i \in \{0,\dots,j-1\}$. By (3.12),
(3.15), the monotonicity of $\phi$, (3.10) and (1.1),
$$\rho(x_{i+1},\bar x) \le \rho(x_{i+1},Tx_i)+\rho(Tx_i,T\bar
x)\le \delta +\phi(\rho(x_i,\bar x))\rho(x_i,\bar x)
$$
$$\le \delta+\phi(\epsilon)\rho(x_i,\bar x)$$
and
$$\rho(x_{i+1},\bar x)-\rho(x_i,\bar x)\le
\delta-(1-\phi(\epsilon))\rho(x_i,\bar x)\le
\delta-(1-\phi(\epsilon))\epsilon\le
-(1-\phi(\epsilon))\epsilon/2.$$ By (3.12) and (3.14) and the
above inequalities
$$-M\le -\rho(x_0,\bar x)\le \rho(x_j,\bar x) -\rho(x_0,\bar x)$$
$$ = \sum_{i=0}^{j-1}[\rho(x_{i+1},\bar x)-\rho(x_i,\bar x)] \le
-j(1-\phi(\epsilon)\epsilon/2) \le
-k(1-\phi(\epsilon))\epsilon/2.$$ This contradicts (3.11). The
contradiction we have reached proves (1.2) and assertion (B).

\head 4. Proof of Theorem 2 \endhead

In the first part of the proof we assume that $T$ is nonexpansive,
i. e., it satisfies (1.1) with $\phi$ identically equal to one.

Let $S\subset [0,1)$ be the following set:
$$S=\{t\in [0,1): tT \text { has a unique fixed
point }x_t\in \text { Int}(G)\}.$$ Since $tT$ is a strict
contraction for each $t\in [0,1)$, it has at most one fixed point.
In order to prove the first part of this theorem, we have to show
that $S=[0,1)$. Since $0\in S$ by assumption and since $[0,1)$ is
connected, it is enough to show that $S$ is both open and closed.

1. $S$ is open: Let $t_0\in S$. From the definition of $S$ it is
clear that $t_0<1$, so there is a real number $q$ such that
$t_0<q<1$. Let $x_{t_0} \in $ Int$(G)$ be the unique fixed point
of $t_0T$.

Since Int$(G)$ is open, there is $r>0$ such that the closed ball
$B[x_{t_0},r]$ of radius $r$ and center $x_{t_0}$ is contained in
Int$(G)$. We have, for all $x\in B[x_{t_0},r]$ and $t \in [0,1)$,
$$||tTx-x_{t_0}|| \le ||tTx-tTx_{t_0}||+|t-t_0|
||Tx_{t_0}||+||t_0Tx_{t_0}-x_{t_0}||$$ $$ \le
t||x-x_{t_0}||+|t-t_0|||Tx_{t_0}|| \le tr+|t-t_0|(||Tx_{t_0}||+1).
\tag4.1$$ Suppose that $t\in [0,1)$ satisfies
$$|t-t_0|<\min\{\frac{r(1-q)}{1+||Tx_{t_0}||},q-t_0\}.\tag4.2$$
 Then $t<q$ and
 $$|t-t_0|\le
\frac{r(1-t)}{1+||Tx_{t_0}||},$$ so $||tTx-x_{t_0}|| \le r$ by
(4.1). Consequently, the closed ball $B[x_{t_0},r]$ is invariant
under $tT$, and
 the Banach fixed point theorem ensures that $tT$ has a unique fixed
point $x_t\in B[x_{t_0},r]\subset$ Int$(G)$. Thus $t\in S$ for all
$t\in [0,1)$ satisfying (4.2).

 2. $S$ is closed: Suppose
$t_0\in [0,1)$ is a limit point of $S$. We have to prove that
$t_0\in S$, and since $0\in S$ we can assume that $t_0>0$. There
is a sequence $(t_n)_n$ in $[0,1)$ such that $t_0=\lim _{n\to
\infty}t_n$, and since $t_0<1$, there is $0<q<1$ such that $t_n<q$
for $n$ large enough. Define
$$ A_0:=\{x_t: t\in
S\cap [0,q]\}.$$ The set $A_0$ is not empty since $0\in A_0$. In
addition, if $t\in S\cap[0,q]$, then
$$
||x_t||=||tTx_t|| \le q(||Tx_t-T0||+||T0||) \le q \phi
(||x_t-0||)||x_t-0||+q||T0||.$$  Therefore
$$
||x_t|| \le  \frac{q||T0||}{1-\phi(||x_t||)q}\le
\frac{||T0||}{1-q}, \tag4.3$$ so $A_0$ is a bounded set, and since
$T$ is Lipschitz, $T(A_0)$ is also bounded, say by $M$. We will
show that $(x_{t_n})_n$ is a Cauchy sequence which converges to
the fixed point $x_{t_0}$ of $t_0T$. Indeed, since $x_{t_n}$ and
$x_{t_m}$ are the fixed points of $t_nT$ and $t_mT$, respectively,
it follows that
$$||x_{t_n}-x_{t_m}||=||t_nTx_{t_n}-t_mTx_{t_m}|| \le
|t_n-t_m|||Tx_{t_n}||+||t_mTx_{t_n}-t_mTx_{t_m}||$$ $$ \le
|t_n-t_m|M+t_m \phi(||x_{t_n}-x_{t_m}||)||x_{t_n}-x_{t_m}||.$$
Hence
$$
||x_{t_n}-x_{t_m}|| \le \frac{|t_n-t_m|M}{1-t_m \phi
(||x_{t_n}-x_{t_m}||)}\le \frac{|t_n-t_m|M}{1-q}. \tag4.4$$

Since $t_n\to t_0$ as $n \to \infty$, we see that $(x_{t_n})_n$ is
indeed Cauchy and hence converges to $x_{t_0}\in \overline{G}$.
Using again the equality $t_nTx_{t_n}=x_{t_n}$, we obtain
$$
||t_0Tx_{t_0}-x_{t_0}|| \le
||t_0Tx_{t_0}-t_0Tx_{t_n}||+||t_0Tx_{t_n}-t_nTx_{t_n}||+
||t_nTx_{t_n}-x_{t_0}||$$
$$=t_0||Tx_{t_0}-Tx_{t_n}||+|t_0-t_n|||Tx_{t_n}||+||x_{t_n}-x_{t_0}||$$
$$\le ||x_{t_0}-x_{t_n}||+|t_0-t_n|M+||x_{t_n}-x_{t_0}||\to
0,$$  so $t_0Tx_{t_0}=x_{t_0}$, i. e., $x_{t_0}$ is indeed a fixed
point of $t_0T$. It remains to show that $x_{t_0}\in$ Int$(G)$,
and this follows from the (L-S) condition, since
 $Tx_{t_0}=\frac{1}{t_0}x_{t_0}$, so (L-S) implies that
 $x_{t_0}\not \in \partial G$ (recall that
$0<t_0<1$). Hence $S$ is closed, as claimed.

The fact that the mapping $t\to x_t$ is Lipschitz on the interval
$[0,b]$ for any $0<b<1$ follows from (4.4).

Suppose now that $T$ satisfies (1.1) with $\phi(t)<1$ for all
positive $t$. Let $(t_n)_n$ be a sequence in $[0,1)$ such that
$t_n\to t_0=1$. The set $A_0$ (and hence the set $T(A_0)$) remain
bounded also when $q=1$, because if $||x_t||\ge 1$, then in (4.3)
we get $||x_t||\le \frac{||T0||}{1-\phi(1)}$, so in any case
$||x_t||\le \max{(1,\frac{||T0||}{1-\phi(1)})}$ (recall that
$\phi(t)<1$). Now, in order to prove that $x_1:=\lim_{t\to
1^{-1}}x_t$ exists, note first that $(x_{t_n})_n$ is Cauchy if
$t_n\to 1$, because otherwise there is $\epsilon>0$ and a
subsequence (call it again $t_n$) such that
$||x_{t_{2n+1}}-x_{t_{2n+2}}||\ge \epsilon$, but from (4.4) we
obtain
$${||x_{t_{2n+1}}-x_{t_{2n+2}}||\le
\frac{|t_{2n+1}-t_{2n+2}|M}{1-t_{2n+2}\phi(\epsilon)}\to 0},$$ a
contradiction. Now, all these  sequences approach the same limit
because for any two such sequences $(x_{t_n})_n,(x_{s_n})_n$, the
interlacing sequence $(t_1,s_1,t_2,s_2\ldots)\to 1$, so
$(x_{t_1},x_{s_1},x_{t_2},x_{s_2},\ldots)$ is also Cauchy. The
fact that $x_1$ is a fixed point of $T$ is proved as above (here,
however,
 one cannot use (L-S) to conclude
that $x_1\in$ Int$(G)$, and indeed it may happen that $x_1\in
\partial G$ as the mapping $T:[-1,\infty)\to R$,
 defined by $Tx=\frac{x-1}{2}$,
 shows).

$$\text { }$$

{\bf Acknowledgments.} The second author was partially supported
by the Fund for the Promotion of Research at the Technion and by
the Technion VPR Fund - B. and G. Greenberg Research Fund
(Ottawa).

\enddocument
\input amstex
\documentstyle {amsppt}
\magnification=\magstephalf \rightheadtext {Metric fixed point
theory} \topmatter
\title Two results in metric fixed point theory
 \endtitle
\author  Daniel Reem, Simeon Reich and Alexander J. Zaslavski \endauthor
\address Department of Mathematics, The Technion-Israel Institute of
Technology, 32000 Haifa, Israel \endaddress \email
dream\@tx.technion.ac.il; sreich\@tx.technion.ac.il; \newline
ajzasl\@tx.technion.ac.il \endemail \abstract We establish two
fixed point theorems for certain mappings of contractive type. The
first result is concerned with the case where such mappings take a
nonempty, closed subset of a complete metric space $X$ into $X$,
and the second with an application of the continuation method to
the case where they satisfy the Leray-Schauder boundary condition
in Banach spaces.
\endabstract
\endtopmatter
\document

\head 1. Introduction  \endhead

In spite of its simplicity (or perhaps because of it), the Banach
fixed point theorem still seems to be the most important result in
metric fixed point theory. As far as we know, the first
significant generalization of Banach's theorem was obtained in
1962 by E. Racotch [5] who replaced Banach's strict contraction
with contractive mappings, that is, with those mappings which
satisfy condition (1.1) below. Such mappings, as well as many
modifications, were studied and used by many authors. Recently, a
renewed interest in contractive mappings has arisen. See, for
example, [6, 7] where genericity and well-posedness results were
established. Another important topic in fixed point theory is the
search for fixed points of nonself-mappings. In the present paper
we combine these two themes by proving two fixed point theorems
for contractive nonself-mappings. In Theorem 1 we present a new
sufficient condition for the existence and approximation of the
unique fixed point of a contractive mapping which maps a nonempty,
closed subset of a complete metric space $X$ into $X$. In Theorem
2 we present an application of the continuation method to
contractive mappings which satisfy the Leray-Schauder boundary
condition on a subset of a Banach space with a nonempty interior.

Let $K$ be a nonempty, closed subset of a complete metric space
$(X,\rho)$. For each $x \in X$ each $r>0$ set
$$B(x,r)=\{y \in X:\; \rho(x,y) \le r\}.$$

\proclaim {Theorem 1}  Assume that $T:K \to X$ satisfies
$$\rho(Tx,Ty) \le \phi(\rho(x,y))\rho(x,y) \text { for each } x,y
\in K, \tag1.1$$
where $\phi:[0,\infty) \to [0,1]$ is a
monotonically decreasing function such that $\phi(t)<1$ for all
$t>0$.

Assume that htere exists a sequence $\{x_n\}_{n=1}^{\infty}
\subset K$ such that $T\bar x=\bar x$.

\endproclaim

\demo {Proof } The uniqueness of $\bar x$ is obviouis. Let us
establish its existence. Let $\epsilon \in (0,1)$. c hoose a
positivw number $\gamma$ such that
$$\gamma <(1-\phi(\epsilon))\epsilon/8.$$
$n_0$ such that $\rho(x_n,Tx_n)<\gamma $ for all integrs $n \ge
n_0$.

Assume that integers $m,n \ge n_0$. We sgow that $\rho(x-m,x_n)
\le \epsiln$. Assume the contrary. Then $\rho(x_m,x_n)> \epsilon$.
Bu

$$\rho(x_n,x_m) \le \rho(x_n,Tx_n)+\rho(Tx_n,Tx_m)+\rho(Tx_m,x_m)
$$
$4\le 2\gama +\phi(\rho(x_n,x_m))\rho(x_m,x_n)\le 3 \gamma
+\phi(\epsilon)\rho(x_m,x_n)$$
$$=\rho(x_m,x_n)-(1-\phi(\epsilon))\rho(x_m,x_n)+2\gamma$$
$$<\rho(x_m,x_n)-(1-\phi(\epsilon))\rho(x_m,x_n)+(1-\phi(\epsilon))\epsilon/4
\le \rho(x_m,x_n)-(1-\phi(\epsilon))\rho(x_m,x_n)/4$$
$$\rho(x_m,x_n)[(3/4)+\phi(\epsilon)/4]<\rho(x_m,x_n),$$
a contradiction.

The contradiction we havereached proves that $\rho(x_m,x_n0 \le
\epsilon$

(P1) For each natural number $n$, there exists $x_n \in K_0$ such
that $T^ix_n$ is defined for all $i=1,\dots,n$.

Then

(A) the mapping $T$ has a unique fixed point $\bar x $ in $K$;

(B) For each $M,\epsilon
>0$, there exist $\delta >0$ and a natural number $k$
such that for each integer $n \ge k$ and each sequence
$\{x_i\}_{i=0}^n \subset K$ satisfying
$$\rho(x_0,\bar x)\le M  \text { and }  \rho(x_{i+1},Tx_i) \le \delta,\;
i=0,\dots,n-1,$$ we have
$$\rho(x_i,\bar x)\le \epsilon, \; i=k,\dots,n.
\tag1.2$$
\endproclaim

Let $G$ be a nonempty subset of a Banach space $(Y,||\cdot||)$. In
[3] J. A. Gatica and W. A. Kirk proved that if $T:\overline{G}\to
Y$ is a strict contraction, then $T$ must have a unique fixed
point $x_1$, under the additional assumptions that the origin is
in the interior Int$(G)$ of $G$ and that $T$ satisfies a certain
boundary condition known as the Leray-Schauder condition:
$$Tx\not = \lambda x \quad \forall x\in
\partial G,\quad \forall\lambda>1. \tag{L-S}$$
Here $G$ is not necessarily convex or bounded. Their proof was
nonconstructive. Later, M. Frigon, A. Granas and Z. E. A. Guennoun
[2], and M. Frigon [1] proved that if $x_t$ is the unique fixed
point of $tT$, then, in fact, the mapping $t\to x_t$ is Lipschitz,
so it gives a partial way to approximate $x_1$. Our second result
extends these theorems to the case where $T$ merely satisfies
(1.1).

\proclaim{Theorem 2} Let $G$ be a nonempty subset of a Banach
space $Y$ with $0\in Int(G)$. Suppose that $T:\overline{G}\to X$
is nonexpansive and that it satisfies condition (L-S). Then for
each $t\in [0,1)$, the mapping $tT:\overline{G}\to X$ has a unique
fixed point $x_t\in Int(G)$ and the mapping $t\to x_t $ is
Lipschitz on $[0,b]$ for any $0<b<1$. If, in addition, $T$
satisfies (1.1), then it has a unique fixed point $x_1\in
\overline{G}$ and the mapping $t\to x_t$ is continuous on $[0,1]$.
In particular, $x_1=\lim_{t\to 1^{-}}x_t$.
\endproclaim

Our paper is organized as follows. Part (A) of Theorem 1 is proved
in Section 2 while Section 3 is devoted to the proof of Theorem
1(B). Finally, we prove Theorem 2 in Section 4.

\head 2. Proof of Theorem 1(A) \endhead

 The uniqueness of $\bar x$ is obvious. To establish its
existence let $x_n \in K_0$ be, for each natural number $n$, the
point provided by property (P1). Fix $\theta_0 \in K$. Since $K_0$
is bounded, there is $c_0>0$ such that
$$\rho(\theta,z) \le c_0 \text { for all } z \in K_0. \tag2.1$$
Let $\epsilon >0$. We will show that there exists a natural number
$k$ such that the following property holds:

(P2) If $n > k$ is an integer and if an integer $i$ satisfies $k
\le i<n$, then
$$\rho(T^ix_n,T^{i+1}x_n) \le \epsilon. \tag2.2$$

Let us assume the contrary. Then for each natural number $k$,
there exist natural numbers $n_k$ and $i_k$ such that
$$k \le i_k<n_k \text { and }
 \rho(T^{i_k}x_{n_k},T^{i_k+1}x_{n_k})>\epsilon.
\tag2.3$$ Choose a natural number $k$ such that
$$k>(\epsilon(1-\phi(\epsilon)))^{-1}(2c_0+\rho(\theta,T\theta)).
\tag2.4$$  By (2.3) and (1.1),
$$\rho(T^ix_{n_k},T^{i+1}x_{n_k})>\epsilon,\; i=0,\dots,i_k.
\tag2.5$$ (Here we use the notation that $T^0z=z$ for al $z \in
K$.)
It follows from (1.1), (2.5) and the monotonicity of $\phi$
that for all $i=0,\dots,i_k-1$,
$$\rho(T^{i+2}x_{n_k},T^{i+1}x_{n_k}) \le
\phi(\rho(T^{i+1}x_{n_k},T^ix_{n_k}))\rho(T^{i+1}x_{n_k},T^ix_{n_k})$$
$$\le \phi(\epsilon)\rho(T^{i+1}x_{n_k},T^ix_{n_k})$$
and
$$\rho(T^{i+2}x_{n_k},T^{i+1}x_{n_k})-\rho(T^{i+1}x_{n_k},T^ix_{n_k})$$
$$\le
(\phi(\epsilon)-1)\rho(T^{i+1}x_{n_k},T^ix_{n_k})<-(1-\phi(\epsilon))\epsilon.
\tag2.6$$ Inequalities (2.6) and (2.3) imply that
$$-\rho(x_{n_k},Tx_{n_k}) \le
\rho(T^{i_k+1}x_{n_k},T^{i_k}x_{n_k})-\rho(x_{n_k},Tx_{n_k})=
$$
$$=\sum_{i=0}^{i_k-1}[\rho(T^{i+2}x_{n_k},T^{i+1}x_{n_k})-
\rho(T^{i+1}x_{n_k},T^ix_{n_k})] \le -(1-\phi(\epsilon)\epsilon)
i_k \le-k(1-\phi(\epsilon)\epsilon
$$
and
$$k(1-\phi(\epsilon))\epsilon \le \rho(x_{n_k},Tx_{n_k}). \tag2.7$$
In view of (2.7), (1.1) and (2.1),
$$k(1-\phi(\epsilon))\epsilon \le \rho(x_{n_k},Tx_{n_k})$$
$$ \le
\rho(x_{n_k},\theta)+\rho(\theta,T\theta)+\rho(T\theta,Tx_{n_k})
\le c_0+\rho(\theta,T\theta)+c_0$$ and
$$k \le
(\epsilon(1-\phi(\epsilon)))^{-1}(2c_0+\rho(\theta,T\theta)).$$
This contradicts (2.4). The contradiction we have reached proves
that for each $\epsilon >0$, there exists a natural number $k$
such that (P2) holds.

Now let $\delta >0$. We show that there exists a  natural number
$k$ such that the following property holds:

(P3) If $n > k$ is an integer and if integers $i,j$ satisfy $k \le
i,j <n$, then
$$\rho(T^ix_n,T^jx_n) \le \delta.$$

To this end, choose a positive number
$$\epsilon <4^{-1}\delta(1-\phi(\delta)). \tag2.8$$
We have already shown that there exists a natural number $k$ such
that (P2) holds.

Assume that natural numbers $n,i$ and $j$ satisfy
$$n > k \text { and } k\le i,j<n. \tag2.9$$
We claim that $\rho(T^ix_n,T^jx_n) \le \delta.$

Assume the contrary. Then
$$\rho(T^ix_n,T^jx_n)>\delta. \tag2.10$$
By (P2), (2.9), (1.1), (2.10) and the monotonicity of $\phi$,
$$\rho(T^ix_n,T^jx_n) \le
\rho(T^ix_n,T^{i+1}x_n)+\rho(T^{i+1}x_n,T^{j+1}x_n)+
\rho(T^{j+1}x_n,T^jx_n)$$ $$ \le
\epsilon+\rho(T^{i+1}x_n,T^{j+1}x_n)+\epsilon \le
2\epsilon+\phi(\rho(T^ix_n,T^jx_n))\rho(T^ix_n,T^jx_n)$$ $$ \le
2\epsilon +\phi(\delta)\rho(T^ix_n,T^jx_n).$$
Together with (2.8)
this implies that
$$\rho(T^ix_n,T^jx_n) \le 2\epsilon(1-\phi(\delta))^{-1}<\delta,$$
a contradiction. Thus we have shown that for each $\delta
>0$, there exists a natural number $k$ such that (P3) holds.

Let $\epsilon >0$. We will show that there exists a natural number
$k$ such that the following property holds:

(P4) If the integers $n_1,n_2 \ge k$, then
$\rho(T^kx_{n_1},T^kx_{n_2}) \le \epsilon$.

Choose a natural number $k$ such that
$$k>((1-\phi(\epsilon))(\epsilon))^{-1}4c_0 \tag2.11$$
and assume that the integers $n_1$ and $n_2$ satisfy $$n_1,n_2 \ge
k. \tag2.12$$ We claim that $\rho(T^kx_{n_1},T^kx_{n_2}) \le
\epsilon$. Assume the contrary. Then
$$\rho(T^kx_{n_1},T^kx_{n_2})>\epsilon.$$
Together with (1.1) this implies that
$$\rho(T^ix_{n_1},T^ix_{n_2})>\epsilon,\; i=0,\dots,k. \tag2.13$$
By (1.1), (2.13) and the monotonicity of $\phi$ we have for
$i=0,\dots,k-1$,
$$\rho(T^{i+1}x_{n_1},T^{i+1}x_{n_2}) \le
\phi(\rho(T^ix_{n_1},T^ix_{n_2}))\rho(T^ix_{n_1},T^ix_{n_2})$$
$$\le \phi(\epsilon)\rho(T^ix_{n_1},T^ix_{n_2})$$
and
$$\rho(T^{i+1}x_{n_1},T^{i+1}x_{n_2})-\rho(T^ix_{n_1},T^ix_{n_2})
\le (\phi(\epsilon)-1)\rho(T^ix_{n_1},T^ix_{n_2}) \le
-(1-\phi(\epsilon))\epsilon.$$ This implies that
$$-\rho(x_{n_1},x_{n_2}) \le
\rho(T^kx_{n_1},T^kx_{n_2})-\rho(x_{n_1},x_{n_2})$$
$$=\sum_{i=0}^{k-1}[\rho(T^{i+1}x_{n_1},T^{i+1}x_{n_2})-\rho(T^ix_{n_1},T^ix_{n_2})]
\le -k(1-\phi(\epsilon))\epsilon.$$ Together with (2.1) this
implies that
$$k(1-\phi(\epsilon))\epsilon \le \rho(x_{n_1},x_{n_2}) \le
\rho(x_{n_1},\theta)+\rho(\theta,x_{n_2})\le 2c_0.$$ This
contradicts (2.11). Thus we have shown that
$$\rho(T^kx_{n_1},T^kx_{n_2})  \le \epsilon.$$
In other words, there exists a natural number $k$ for which (P4)
holds.

Let $\epsilon >0$. By (P4), there exists a natural number $k_1$
such that
$$\rho(T^{k_1}x_{n_1},T^{k_1}x_{n_2}) \le \epsilon/4 \text  { for
all integers } n_1,n_2 \ge k_1. \tag2.14$$ By (P3), there exists a
natural number $k_2$ such that $$\rho(T^ix_n,T^jx_n) \le
\epsilon/4 \text { for all natural numbers } n,j,i \text {
satisfying } k_2 \le i,j< n. \tag2.15$$

Assume now that the natural numbers $n_1,n_2, i$ and $j$ satisfy
$$n_1,n_2 > k_1+k_2,\; i,j \ge k_1+k_2,\; i<n_1,\; j <
n_2. \tag2.16$$ We claim that
$$\rho(T^ix_{n_1},T^jx_{n_2}) \le
\epsilon.$$
By (2.14), (2.16) and (1.1)
$$\rho(T^{k_1+k_2}x_{n_1},T^{k_1+k_2}x_{n_2}) \le
\rho(T^{k_1}x_{n_1},T^{k_1}x_{n_2}) \le \epsilon/4. \tag2.17$$
In
view of (2.16) and (2.15),
$$\rho(T^{k_1+k_2}x_{n_1},T^ix_{n_1}) \le \epsilon/4 \text { and }
\rho(T^{k_1+k_2}x_{n_2},T^jx_{n_2})\le \epsilon/4.$$ Together with
(2.17) these inequalities imply that
$$\rho(T^ix_{n_1},T^jx_{n_2}) \le
\rho(T^ix_{n_1},T^{k_1+k_2}x_{n_1})$$ $$+
\rho(T^{k_1+k_2}x_{n_1},T^{k_1+k_2}x_{n_2})+
\rho(T^{k_1+k_2}x_{n_2},T^jx_{n_2}) <\epsilon.$$ Thus we have
shown that the following property holds:

(P5) For each $\epsilon
>0$, there exists a natural number $k(\epsilon)$ such that
$$\rho(T^ix_{n_1},T^jx_{n_2})\le \epsilon $$ for all natural numbers
$n_1,n_2 \ge k(\epsilon)$, $i \in [k(\epsilon),n_1)$ and $j \in
[k(\epsilon),n_2).$

Consider the two sequences $\{T^{n-2}x_n\}_{n=2}^{\infty}$ and
$\{T^{n-1}x_n\}_{n=2}^{\infty}$. Property (P5) implies that both
of them are Cauchy and that
$$\lim_{n
\to\infty}\rho(T^{n-1}x_n,T^{n-2}x_n)=0.$$ Therefore there exists
$\bar x \in K$ such that
$$\lim_{n \to \infty} \rho(\bar x,T^{n-2}x_n)=\lim_{n \to
\infty}\rho(\bar x,T^{n-1}x_n)=0.$$ Since the mapping $T$ is
continuous, $T\bar x=\bar x$ and assertion (A) is proved.

\head 3. Proof of Theorem 1(B) \endhead

For each $x \in X$ and $r>0$, set
$$B(x,r)=\{y \in X:\; \rho(x,y) \le r\}. \tag3.1$$
Choose $\delta_0>0$ such that
$$\delta_0<M(1-\phi(M/2))/4.\tag3.2$$
Assume that
$$y \in K\cap B(\bar x,M),\; z \in X \text { and } \rho(z,Ty)\le \delta_0.
\tag3.3$$ By (3.3) and (1.1),
$$\rho(\bar x,z)\le \rho(\bar x,Ty)+\rho(Ty,z)\le \rho(T\bar
x,Ty)+\delta_0\le \phi(\rho(\bar x,y))\rho(\bar x,y)+\delta_0.
\tag3.4$$ There are two cases:
$$\rho(y,\bar x) \le M/2; \tag3.5$$
$$\rho(y,\bar x)>M/2. \tag3.6$$
Assume that (3.5) holds. By (3.4), (3.5) and (3.2),
$$\rho(\bar x,z) \le \rho(\bar x,y)+\delta_0\le M/2+\delta_0<M.
\tag3.7$$ If (3.6) holds, then by (3.4), (3.3), (3.6) and the
monotonicity of $\phi$,
$$\rho(\bar x,z)\le \delta_0+\phi(M/2)\rho(\bar x,y)\le
\delta_0+\phi(M/2)M$$ $$<(M/4)(1-\phi(M/2))+\phi(M/2)M\le M.$$
Thus $\rho(\bar x,z) \le M$ in both cases.

We have shown that
$$\rho(\bar x,z)\le  M \text { for each } z \in X \text { and } y
\in K \cap B(\bar x,M) \text { satisfying } \rho(z,Ty)\le
\delta_0.\tag3.8$$ Since $M$ is any positive number, we conclude
that there is $\delta_1>0$ such that
$$\rho(\bar x,z)\le \epsilon \text { for each } z \in X \text {
and } y \in K \cap B(\bar x,\epsilon) \text { satisfying }
\rho(z,Ty)\le \delta_1. \tag3.9$$ Choose a positive number $\delta
$ such that
$$\delta <\min\{\delta_0,\delta_1,\;
\epsilon(1-\phi(\epsilon))4^{-1}\} \tag3.10$$ and a natural number
$k$ such that
$$k>4(M+1)(1-\phi(\epsilon)\epsilon)^{-1}+4. \tag3.11$$
Let $n \ge k$ be a natural number and assume that $\{x_i\}_{i=0}^n
\subset K$ satisfies
$$\rho(x_0,\bar x) \le M  \text { and } \rho(x_{i+1},Tx_i) \le
\delta,\;i=0,\dots,n-1. \tag3.12$$ We claim that (1.2) holds. By
(3.8), (3.12) and the inequality $\delta <\delta_0$ (see (3.10)),
$$\{x_i\}_{i=0}^k \subset B(\bar x,M). \tag3.13$$
Assume that (1.2) does not hold. Then there is an integer $j$ such
that
$$j \in \{k,n\} \text { and } \rho(x_j,\bar x)>\epsilon. \tag3.14$$
By (3.14), (3.12), (3.9) and (3.10),
$$\rho(x_i,\bar x)>\epsilon,\;
i=0,\dots,j. \tag3.15$$ Let $i \in \{0,\dots,j-1\}$. By (3.12),
(3.15), the monotonicity of $\phi$, (3.10) and (1.1),
$$\rho(x_{i+1},\bar x) \le \rho(x_{i+1},Tx_i)+\rho(Tx_i,T\bar
x)\le \delta +\phi(\rho(x_i,\bar x))\rho(x_i,\bar x)
$$
$$\le \delta+\phi(\epsilon)\rho(x_i,\bar x)$$
and
$$\rho(x_{i+1},\bar x)-\rho(x_i,\bar x)\le
\delta-(1-\phi(\epsilon))\rho(x_i,\bar x)\le
\delta-(1-\phi(\epsilon))\epsilon\le
-(1-\phi(\epsilon))\epsilon/2.$$ By (3.12) and (3.14) and the
above inequalities
$$-M\le -\rho(x_0,\bar x)\le \rho(x_j,\bar x) -\rho(x_0,\bar x)$$
$$ = \sum_{i=0}^{j-1}[\rho(x_{i+1},\bar x)-\rho(x_i,\bar x)] \le
-j(1-\phi(\epsilon)\epsilon/2) \le
-k(1-\phi(\epsilon))\epsilon/2.$$ This contradicts (3.11). The
contradiction we have reached proves (1.2) and assertion (B).

\head 4. Proof of Theorem 2 \endhead

In the first part of the proof we assume that $T$ is nonexpansive,
i. e., it satisfies (1.1) with $\phi$ identically equal to one.

Let $S\subset [0,1)$ be the following set:
$$S=\{t\in [0,1): tT \text { has a unique fixed
point }x_t\in \text { Int}(G)\}.$$
Since $tT$ is a strict
contraction for each $t\in [0,1)$, it has at most one fixed point.
In order to prove the first part of this theorem, we have to show
that $S=[0,1)$. Since $0\in S$ by assumption and since $[0,1)$ is
connected, it is enough to show that $S$ is both open and closed.

1. $S$ is open: Let $t_0\in S$. From the definition of $S$ it is
clear that $t_0<1$, so there is a real number $q$ such that
$t_0<q<1$. Let $x_{t_0} \in $ Int$(G)$ be the unique fixed point
of $t_0T$.

Since Int$(G)$ is open, there is $r>0$ such that the closed ball
$B[x_{t_0},r]$ of radius $r$ and center $x_{t_0}$ is contained in
Int$(G)$. We have, for all $x\in B[x_{t_0},r]$ and $t \in [0,1)$,
$$||tTx-x_{t_0}|| \le ||tTx-tTx_{t_0}||+|t-t_0|
||Tx_{t_0}||+||t_0Tx_{t_0}-x_{t_0}||$$ $$ \le
t||x-x_{t_0}||+|t-t_0|||Tx_{t_0}|| \le tr+|t-t_0|(||Tx_{t_0}||+1).
\tag4.1$$ Suppose that $t\in [0,1)$ satisfies
$$|t-t_0|<\min\{\frac{r(1-q)}{1+||Tx_{t_0}||},q-t_0\}.\tag4.2$$
 Then $t<q$ and
 $$|t-t_0|\le
\frac{r(1-t)}{1+||Tx_{t_0}||},$$ so $||tTx-x_{t_0}|| \le r$ by
(4.1). Consequently, the closed ball $B[x_{t_0},r]$ is invariant
under $tT$, and
 the Banach fixed point theorem ensures that $tT$ has a unique fixed
point $x_t\in B[x_{t_0},r]\subset$ Int$(G)$. Thus $t\in S$ for all
$t\in [0,1)$ satisfying (4.2).

 2. $S$ is closed: Suppose
$t_0\in [0,1)$ is a limit point of $S$.
We have to prove that
$t_0\in S$, and since $0\in S$ we can assume that $t_0>0$. There
is a sequence $(t_n)_n$ in $[0,1)$ such that $t_0=\lim _{n\to
\infty}t_n$, and since $t_0<1$, there is $0<q<1$ such that $t_n<q$
for $n$ large enough. Define
$$ A_0:=\{x_t: t\in
S\cap [0,q]\}.$$ The set $A_0$ is not empty since $0\in A_0$. In
addition, if $t\in S\cap[0,q]$, then
$$
||x_t||=||tTx_t|| \le q(||Tx_t-T0||+||T0||) \le q \phi
(||x_t-0||)||x_t-0||+q||T0||.$$  Therefore
$$
||x_t|| \le  \frac{q||T0||}{1-\phi(||x_t||)q}\le
\frac{||T0||}{1-q}, \tag4.3$$ so $A_0$ is a bounded set, and since
$T$ is Lipschitz, $T(A_0)$ is also bounded, say by $M$. We will
show that $(x_{t_n})_n$ is a Cauchy sequence which converges to
the fixed point $x_{t_0}$ of $t_0T$. Indeed, since $x_{t_n}$ and
$x_{t_m}$ are the fixed points of $t_nT$ and $t_mT$, respectively,
it follows that
$$||x_{t_n}-x_{t_m}||=||t_nTx_{t_n}-t_mTx_{t_m}|| \le
|t_n-t_m|||Tx_{t_n}||+||t_mTx_{t_n}-t_mTx_{t_m}||$$ $$ \le
|t_n-t_m|M+t_m \phi(||x_{t_n}-x_{t_m}||)||x_{t_n}-x_{t_m}||.$$
Hence
$$
||x_{t_n}-x_{t_m}|| \le \frac{|t_n-t_m|M}{1-t_m \phi
(||x_{t_n}-x_{t_m}||)}\le \frac{|t_n-t_m|M}{1-q}. \tag4.4$$

Since $t_n\to t_0$ as $n \to \infty$, we see that $(x_{t_n})_n$ is
indeed Cauchy and hence converges to $x_{t_0}\in \overline{G}$.
Using again the equality $t_nTx_{t_n}=x_{t_n}$, we obtain
$$
||t_0Tx_{t_0}-x_{t_0}|| \le
||t_0Tx_{t_0}-t_0Tx_{t_n}||+||t_0Tx_{t_n}-t_nTx_{t_n}||+
||t_nTx_{t_n}-x_{t_0}||$$
$$=t_0||Tx_{t_0}-Tx_{t_n}||+|t_0-t_n|||Tx_{t_n}||+||x_{t_n}-x_{t_0}||$$
$$\le ||x_{t_0}-x_{t_n}||+|t_0-t_n|M+||x_{t_n}-x_{t_0}||\to
0,$$  so $t_0Tx_{t_0}=x_{t_0}$, i. e., $x_{t_0}$ is indeed a fixed
point of $t_0T$. It remains to show that $x_{t_0}\in$ Int$(G)$,
and this follows from the (L-S) condition, since
 $Tx_{t_0}=\frac{1}{t_0}x_{t_0}$, so (L-S) implies that
 $x_{t_0}\not \in \partial G$ (recall that
$0<t_0<1$). Hence $S$ is closed, as claimed.

The fact that the mapping $t\to x_t$ is Lipschitz on the interval
$[0,b]$ for any $0<b<1$ follows from (4.4).

Suppose now that $T$ satisfies (1.1) with $\phi(t)<1$ for all
positive $t$. Let $(t_n)_n$ be a sequence in $[0,1)$ such that
$t_n\to t_0=1$. The set $A_0$ (and hence the set $T(A_0)$) remain
bounded also when $q=1$, because if $||x_t||\ge 1$, then in (4.3)
we get $||x_t||\le \frac{||T0||}{1-\phi(1)}$, so in any case
$||x_t||\le \max{(1,\frac{||T0||}{1-\phi(1)})}$ (recall that
$\phi(t)<1$). Now, in order to prove that $x_1:=\lim_{t\to
1^{-1}}x_t$ exists, note first that $(x_{t_n})_n$ is Cauchy if
$t_n\to 1$, because otherwise there is $\epsilon>0$ and a
subsequence (call it again $t_n$) such that
$||x_{t_{2n+1}}-x_{t_{2n+2}}||\ge \epsilon$, but from (4.4) we
obtain
$${||x_{t_{2n+1}}-x_{t_{2n+2}}||\le
\frac{|t_{2n+1}-t_{2n+2}|M}{1-t_{2n+2}\phi(\epsilon)}\to 0},$$ a
contradiction. Now, all these  sequences approach the same limit
because for any two such sequences $(x_{t_n})_n,(x_{s_n})_n$, the
interlacing sequence $(t_1,s_1,t_2,s_2\ldots)\to 1$, so
$(x_{t_1},x_{s_1},x_{t_2},x_{s_2},\ldots)$ is also Cauchy. The
fact that $x_1$ is a fixed point of $T$ is proved as above (here,
however,
 one cannot use (L-S) to conclude
that $x_1\in$ Int$(G)$, and indeed it may happen that $x_1\in
\partial G$ as the mapping $T:[-1,\infty)\to R$,
 defined by $Tx=\frac{x-1}{2}$,
 shows).

$$\text { }$$

{\bf Acknowledgments.} The second author was partially supported
by the Fund for the Promotion of Research at the Technion and by
the Technion VPR Fund - B. and G. Greenberg Research Fund
(Ottawa).

\enddocument